\documentclass[12pt]{amsart}
\usepackage{a4wide}
\usepackage{amsmath,amssymb,amsfonts,amsthm,dutchcal}
\usepackage{tikz}
\usetikzlibrary{decorations.pathreplacing}
\usepackage{xcolor}
\usepackage{graphicx}

\usepackage{cancel}

\usepackage[colorlinks=true,urlcolor=blue,
citecolor=red,linkcolor=blue,linktocpage,pdfpagelabels,
bookmarksnumbered,bookmarksopen]{hyperref}

\theoremstyle{plain}

\hypersetup{
 colorlinks   = true,
 urlcolor     = blue,
 linkcolor    = blue,
 citecolor   = red ,
 bookmarksopen=true
}
\usepackage[noabbrev]{cleveref}

\theoremstyle{plain}
\newtheorem{thrm}{Theorem}[section]
\newtheorem{lemma}[thrm]{Lemma}

\newtheorem{cor}[thrm]{Corollary}

\newtheorem{dfn}[thrm]{Definition}

\allowdisplaybreaks
\begin{document}
\newcommand{\psr}{P^-}
\newcommand{\qr}{Q_\rho}
\newcommand{\sn}{\mathbb{S}^{n-1}}
\newcommand{\SL}{\mathcal L^{1,p}( D)}
\newcommand{\Lp}{L^p( Dega)}
\newcommand{\CO}{C^\infty_0( \Omega)}
\newcommand{\Rn}{\mathbb R^n}
\newcommand{\Rm}{\mathbb R^m}
\newcommand{\R}{\mathbb R}
\newcommand{\Om}{\Omega}
\newcommand{\Hn}{\mathbb H^n}
\newcommand{\aB}{\alpha B}

\newcommand{\BVX}{BV_X(\Omega)}
\newcommand{\p}{\partial}
\newcommand{\I}{\mathbb I}
\newcommand{\bG}{\boldsymbol{G}}
\newcommand{\bg}{\mathfrak g}
\newcommand{\bz}{\mathfrak z}
\newcommand{\bv}{\mathfrak v}
\newcommand{\Bux}{\mbox{Box}}
\newcommand{\X}{\mathcal X}
\newcommand{\Y}{\mathcal Y}
\newcommand{\B}{\mathcal{B}}
\newcommand{\K}{\mathcal{K}}
\newcommand{\osc}{\operatorname{osc}}

\numberwithin{equation}{section}

\title[Hölder Continuity etc.]{Hölder Continuity and Harnack  estimate for  non-homogeneous  parabolic equations}

\author{Vedansh Arya}
\address{Department of Mathematics and Statistics\\
University of Jyväskylä\\ Finland}\email[Vedansh Arya]{vedansh.v.arya@jyu.fi}

\author{Vesa Julin}
\address{Department of Mathematics and Statistics\\
University of Jyväskylä\\ Finland}\email[Vesa Julin]{vesa.julin@jyu.fi}

%

%
%
%

\begin{abstract}

In this paper we continue the study on intrinsic Harnack inequality for  non-homogeneous parabolic equations in non-divergence form initiated by the first author in \cite{A}. 
We establish a forward-in-time intrinsic Harnack inequality, which in particular implies the H\"{o}lder continuity of the solutions. We also provide a Harnack type estimate on global scale which quantifies the strong minimum  principle. In the time-independent setting, this together with  \cite{A} provides an alternative proof of the generalized Harnack inequality proven by the second author in \cite{j}.  
\end{abstract}
\maketitle

 
\section{Introduction and the statement of the main results}
In this paper we consider parabolic equation of the type
\begin{align}\label{maineq}
    F(D^2u, Du, x, t)-u_t=0,
\end{align}
where  $F$ is  uniformly elliptic w.r.t. the Hessian and has  a nonlinear growth w.r.t. the gradient. More precisely,  we assume that  there exist constants $0 < \lambda \leq \Lambda$ such that 
\begin{align*}
    \lambda Tr(N) \leq F(M+N,p,x,t)-F(M,p,x,t) \leq \Lambda Tr(N)
\end{align*}
for all symmetric matrices $M,N$ with $  N \geq 0$  and for every $(p,x,t) \in \R^n \times Q_r(x_0,t_0)$. Here  $Q_{r}(x_0,t_0)$ is a cube centered at  $(x_0, t_0)\in \R^{n+1}$.  We assume further that $F$ has the following growth in the gradient variable,
\begin{align}\label{nonlin}
    |F(0,p,x,t)| \leq \phi(|p|)
\end{align}
for every $(p,x,t) \in \R^n \times Q_r(x_0,t_0)$, where  $\phi : [0,\infty) \rightarrow [0,\infty)$ is of the form $\phi(t)=\eta(t)t$ and satisfies the following conditions:
\begin{enumerate} 
    \item[(P1)] $\phi : [0,\infty) \rightarrow [0,\infty)$ is increasing, locally Lipschitz continuous in $(0,\infty)$ and $\phi(t)\geq t$ for every $t\geq0$. Moreover, $\eta : (0,\infty) \rightarrow [1,\infty)$ is nonincreasing on $(0,1)$ and nondecreasing on $[1,\infty)$;
    \item[(P2)] $\eta$ satisfies 
    \begin{equation} \label{eq:P2}
        \underset{t\rightarrow\infty}{\text{lim}}\frac{t\eta^{'}(t)}{\eta(t)}\text{log}(\eta(t))=0;
    \end{equation}
    \item[(P3)] There is a constant $\Lambda_0$ such that 
    \begin{equation} \label{eq:P3}
        \eta(st)\leq \Lambda_0\eta(s)\eta(t);
     \end{equation}
    for every $s$, $t\in (0,\infty)$.
\end{enumerate}

The elliptic case is studied by the second author in \cite{j}, and there the main problem is to find a Harnack type estimate which quantifies the  strong minimum principle whenever it is true  and generalizes the Harnack inequality from the homogeneous setting \cite{Ca,KS}. The related boundary problem  is studied in \cite{AJ, Lu, LOT}. In this paper we first discuss the assumptions (P1)-(P3).  In particular, the condition (P2) in \eqref{eq:P2} might first seem rather resctrictive. We note that (P2) roughly states that the function $\eta$ is \emph{slowly varying} \footnote{A positive function $\eta$ is slowly varying if for all $\lambda > 0$ it holds $\lim_{s \to \infty} \frac{\eta(\lambda s)}{\eta(s)} = 1$.}  in a quantitative way  and, in particular, it implies the slow growth estimate on $\eta$, i.e., for every $\varepsilon>0$ there is $C_\varepsilon$ such that $\eta(t) \leq C_\varepsilon t^{\varepsilon}$ for $t \geq 0$ \cite[Proposition 1.3.6]{BGT}. Our first result (Theorem \ref{thm1}) justifies the assumption  (P2) in the sense that we show that the Harnack estimate in \cite{j} is false if the non-linearity is of  the form  $\eta(t) = \max \{t^\varepsilon, 1\}$ for any $\varepsilon>0$.

The parabolic case is studied by the first author in \cite{A}, where the main result is  the backward-in-time intrinsic Harnack inequality. The parabolic case is the main focus of this paper and we first provide the analogous  forward-in-time intrinsic Harnack inequality (Theorem \ref{fullharnack}), and show then that this intrinsic Harnack inequality implies  the H\"older continuity of the solutions (Theorem \ref{Th:holder}). The problem has some similarities to the p-parabolic equation which is also non-nomogeneous, and where the concept of intrinsic Harnack inequality is developed \cite{DGV1, DGV3} (see also \cite{Ku}). The difference is that in the p-parabolic equation the elliptic operator has a different scaling  than the time derivative term, while in  \eqref{maineq} it is the elliptic operator that does not have any homogeneity properties. Therefore, from the point of view of the scaling,  the equation  \eqref{maineq} is  more similar to the reaction diffusion equation than to the p-parabolic. However, to the best our knowledge  these results are not known for equations of type
\[
u_t = F(D^2u) +f(u).
\]

The intrinsic Harnack inequality is a result on intrinsic scale and it is a relevant question if one can have a generalized Harnack inequality for non-negative solutions of \eqref{maineq} similar to  the elliptic one in \cite{j}. In our last result (Theorem \ref{Th:int_h}) we provide a Harnack type estimate for non-negative solutions of \eqref{maineq} on global scale which quantifies the strong  minimum principle.

\subsection{Statement of the main results}

In order to state our main result we first simplify the setting and, following the ideas of Caffarelli \cite{Ca, cc},   replace the equation \eqref{maineq} by two extremal inequalities which take into account the   ellipticity assumption and the  growth condition of the  drift term. To be more precise,  we assume that $u \in C(Q_2(x_0,t_0))$ is a viscosity supersolution of
\begin{equation}\label{eq1}
    P_{\lambda, \Lambda}^-(D^2u)-u_t\leq\phi(|Du|)
\end{equation}
and  a viscosity subsolution of 
\begin{equation}\label{eq2}
    P_{\lambda, \Lambda}^+(D^2u)-u_t\geq-\phi(|Du|),
\end{equation}
where $P_{\lambda, \Lambda}^{\pm}$  denotes the extremal Pucci operators defined later in \eqref{pucci}. We refer to Section \ref{s:n} for the precise notion of viscosity sub- and supersolutions. The elliptic counterpart to \eqref{eq1} and \eqref{eq2} is defined analogously. 

Throughout the paper we say that a constant $C$ is \emph{universal}  if it  depends only  on the  ellipticity constants,  the nonlinearity $\phi$ and on the  dimension $n$.

We recall that in \cite{j} it is proven that in the elliptic case   if $u \in C(B_2)$  is a non-negative viscosity supersolution of $ P_{\lambda,\Lambda}^-(D^2u) \leq\phi(|Du|)$ and a viscosity subsolution $ P_{\lambda,\Lambda}^+(D^2u) \geq-\phi(|Du|)$ in $B_2$, then there is  a universal constant $C$ such that for $m = \inf_{B_1} u$ and $M = \sup_{B_1}u$  it holds
\begin{equation}\label{eq:harnack1}
   \int_{m}^{M} \frac{ds}{\phi(s)} \leq C. 
\end{equation}
It is not difficult to see that if $u$ is as above, then it satisfies the strong minimum principle when $\phi$ satisfies the so called Osgood condition \cite{j}
\begin{equation}\label{eq:osgood}
   \int_{0}^{1} \frac{ds}{\phi(s)} = \infty. 
\end{equation}
Therefore it is clear that \eqref{eq:harnack1} quantifies the strong minimum principle when \eqref{eq:osgood} holds.

Our first result justifies the assumption (P2) in \eqref{eq:P2}. 
\begin{thrm}
\label{thm1}
Assume $n \geq 3$ and fix any $\varepsilon_0>0$. Then there exist $\lambda,\Lambda >0$ and a sequence of positive functions $u_k \in C^2(B_2)$ such that  
\[
 \big|P_{\lambda, \Lambda}^- (D^2u_k)\big| \leq |Du_k|^{1+\varepsilon_0} \qquad \text{point-wise in } \, B_2,
\]
while $m_k = \inf_{B_1} u_k \to 0$ as $k \to \infty$ and $M_k =  \sup_{B_1} u_k \geq 1$ for all $k \in \mathbb{N}$. In particular, it holds 
\[
\int_{m_k}^{M_k} \frac{ds}{s^{1+\varepsilon_0} +s}  \to \infty \qquad \text{as }\, k \to \infty. 
\]
\end{thrm}
We note that the condition (P2) in \eqref{eq:P2} is an assumption only on the asymptotic behavior of the nonlinearity $\phi$ and plays no role whether the Osgood condition \eqref{eq:osgood} holds or not. The point of Theorem \ref{thm1} is that too weak information on the asymptotic growth of $\phi$ is the reason for the Harnack estimate to fail, even if the Harnack estimate means  to quantify only  the strong minimum principle.

Next we turn our attention to the parabolic setting  \eqref{maineq}, which is the main focus of this paper. We recall the following parabolic Harnack type inequality  proven in \cite{A}  for solutions to \eqref{maineq} in a  suitable intrinsic geometry    corresponding to the nonlinearity.  Let $u \in C(Q_2)$ be a positive viscosity supersolution of \eqref{eq1} and viscosity subsolution of \eqref{eq2}. There is a universal constant $C>0$ such that 
\begin{equation}\label{Harnack}
    \underset{A_1^-}{\text{sup}}\hspace{0.8mm} u(a_0 x, a_0^2 t) \leq  C u(0,0) \hspace{2mm}\text{for} \hspace{2mm} a_0 \le  \frac{u(0,0)}{C(\phi(u(0,0))+u(0,0))},
\end{equation}
where $A_1^-=\Big\{(x,t):|x|_{\infty} \leq \frac{c_n}{2}, -1+ \frac{c_n^2}{4} \leq t \leq -1+\frac{c_n^2}{2}\Big\}$ and $|\cdot|_\infty$ is the $l^\infty$-norm. Here $c_n\leq 1$ depends only on $n$. The estimate \eqref{Harnack} is an intrinsic Harnack estimate as the scaling on the LHS depends on the value of the solution at the base point which is similar to  \cite{DGV3}.

Our next result is the forward-in-time counterpart of \eqref{Harnack} and thus we obtain the  complete intrinsic Harnack estimate for \eqref{maineq}. To this aim we denote 
\begin{equation}\label{def:Aplus}
A^+_{\rho}(x_0,t_0):=\Big\{(x,t):|x-x_0|_{\infty} < \frac{\rho c_n}{2}, \rho^2- \frac{(\rho c_n)^2}{2} < t-t_0 \leq \rho^2-\frac{(\rho c_n)^2}{4}\Big\},
\end{equation}
where $c_n\leq 1$ is as above. 
\begin{thrm}\label{fullharnack}
Let $u \in C(Q_2)$ be a non-negative viscosity supersolution of \eqref{eq1} and viscosity subsolution of \eqref{eq2}. There is a universal constant $C>0$ such that if $u(x_0,t_0) >0$ and $Q_{2\rho}(x_0,t_0) \cup A^+_{\rho}(x_0,t_0) \subset Q_2$ for $$\rho \le \alpha_0:=\frac{u(x_0,t_0)}{C(\phi(u(x_0,t_0))+u(x_0,t_0))}$$ then it holds
\begin{equation}\label{fharnackp}
    u(x_0,t_0) \le C\underset{A_{\rho}^+(x_0,t_0)}{\inf}\hspace{0.8mm} u( x, t).
\end{equation}
\end{thrm}
The proof of Theorem \ref{fullharnack}  relies on  \eqref{Harnack} and a  careful continuity-type argument.

As in the p-parabolic case  \cite{ DGV1, DGV3},  using the estimate  \eqref{fharnackp} we obtain the H\"older continuity of the solutions of \eqref{maineq}. 
\begin{thrm}\label{Th:holder}
Let $u \in C(Q_2)$ be a viscosity supersolution of \eqref{eq1} and viscosity subsolution of \eqref{eq2}. Then, there exists a universal constant $\alpha \in (0,1)$ such that $u$ is locally $\alpha$-H\"older continuous. More precisely, there is a universal  constant $C>1$,  such that for all $(x,t), (y,s) \in Q_1$, we have 
\begin{equation}\label{thm:holder}
|u(x,t)-u(y,s)|  \le C\phi(\|u\|_{L^{\infty}(Q_2)}) (|x-y|+|t-s|)^{\alpha}.
\end{equation}
\end{thrm}
We note that in \eqref{thm:holder} the exponent $\alpha \in (0,1)$ is universal, but the H\"older norm depends on the solution $u$ in a nonlinear way. This is necessary already in the elliptic case \cite{Tr}. Our  proof for Theorem \ref{Th:holder}  draws inspiration from the p-parabolic case \cite{DGV3}. However,  since the elliptic operator in  the equation \eqref{maineq}  does not have any degree of homogeneity, the intrinsic cylinders in Theorem   \ref{fullharnack}  may not be quantitatively monotone in size. This causes  challenges which we need to overcome in  order to obtain a uniform $\alpha>0$.  

Finally, we study whether we may obtain a Harnack type estimate on global scale,  which quantifies the strong minimum principle similar to \eqref{eq:harnack1}. Recall that the strong minimum principle is related to the Osgood condition \eqref{eq:osgood}. 

\begin{thrm}\label{Th:int_h}
Let $u \in C(Q_{4}(0,1))$ be a non-negative viscosity supersolution
of \eqref{eq1} and viscosity subsolution of \eqref{eq2} in $Q_{4}(0,1)$. There exist a universal constant $C>1$ and time levels $t_1,t_2\geq 0 $, depending on the value $u(0,0)$ such that  $\frac{1}{C(\eta(u(0,0))+1)^2} \le t_1,t_2 \le 1$, so that it holds 
\begin{equation} \label{eq:para-harnack}
 \int_{\inf_{B_1}u(\cdot,t_1)}^{\sup_{B_1}u(\cdot,-t_2)} \frac{ds}{ \phi(s)} \leq C.
 \end{equation}
\end{thrm}
The parabolic Harnack estimate \eqref{eq:para-harnack} takes the same form as the elliptic one \eqref{eq:harnack1}, with the difference that here we have the waiting times $t_1,t_2$ as usual with parabolic equations \cite{is, W1}. Since the scaling of the equation \eqref{maineq} does not have any monotonicity properties, we can only give an estimate  for the waiting times. However, the proof of Theorem \ref{Th:int_h} implies that, under the Osgood-condition \eqref{eq:osgood}, if $u(0,0)= 0$ then the  waiting time $t_2$ can  be chosen  zero. We state this in the following corollary, which is  a quantification of the strong minimum principle at the time level $t=0$.
\begin{cor}
\label{cor:harnack}
Let $u \in C(Q_{4})$ be as in Theorem \ref{Th:int_h} and assume that $\phi$ satisfies  \eqref{eq:osgood}. If $u(0,0)= 0$ then 
\[
u(x,0)= 0 \qquad \text{for all }\, x \in B_1.
\]
\end{cor}

In the next section  we introduce some basic notations and in Section \ref{s:proof} we prove our main results.

\section{Notations and Preliminaries}\label{s:n}
We denote a point in space by $x \in \R^n$ and in  space-time  by $(x,t) \in \R^n \times \R$.  We denote the  Euclidean norm of $x$ in $\R^n$ by $|x|$ and the $l^\infty$-norm by $|x|_{\infty} = \text{max}\{ |x_1|, |x_2|,..., |x_n|\}$. 
We denote $Df$  the gradient of $f$ in $x$,  $f_t$  with respect to time $t$ and $D^2f$ denotes the Hessian matrix of $f$ with respect to $x$. We denote the  ball of radius $r$ centered at $x_0$ by $B_r(x_0)$
and its closure by  $\overline{B}_r(x_0).$ 
A cube of radius $\rho$ centered at $(x_0,t_0)$ is defined as 
\begin{align*}
    \qr(x_0,t_0) = \{x \in \R^n:|x-x_0|_{\infty} < \rho \}\times (t_0-\rho^2,t_0].
\end{align*}
We denote $B_r$ and $Q_r$ if the ball and the cube are centered at the origin. 

Let $S$ be the space of real $n\times n$ symmetric matrices. We recall the definition of Pucci's extremal operators (for more detail see \cite{cc}). For $M\in S$, Pucci's extremal operators with ellipticity constant $0<\lambda\leq\Lambda$ are defined as
\begin{align}\label{pucci}
    P_{\lambda,\Lambda}^-(M)=\lambda\sum_{e_i>0}e_i+\Lambda\sum_{e_i<0}e_i,
    \\
    P_{\lambda,\Lambda}^+(M)=\Lambda\sum_{e_i>0}e_i+\lambda\sum_{e_i<0}e_i\notag,
\end{align}
where $e_i$'s are the eigenvalues of $M$.

We recall the definition of a viscosity supersolution of \eqref{eq1} and  viscosity subsolution of \eqref{eq2}.
\begin{dfn}
A lower semicontinuous function $u :Q_{r}(x_0,t_0) \rightarrow \R$ is a viscosity supersolution  of \eqref{eq1}  in $Q_r(x_0,t_0)$ if the following holds: if $(x,t)\in Q_r(x_0,t_0)$ and $\varphi \in C^2(Q_r(x_0,t_0))$ are such that $\varphi\leq u$ and $\varphi(x,t)=u(x,t)$ then
\begin{align*}
    P_{\lambda,\Lambda}^-(D^2\varphi(x,t))-\varphi_t(x,t)\leq \phi(|D \varphi(x,t)|).
\end{align*}
An upper semicontinuous function $u :Q_r(x_0,t_0) \rightarrow \R$ is a viscosity subsolution  of \eqref{eq2}  in $Q_r(x_0,t_0)$ if the following holds: if $(x,t)\in Q_r(x_0,t_0)$ and $\varphi \in C^2(Q_r(x_0,t_0))$ are such that $\varphi\geq u$ and $\varphi(x,t)=u(x,t)$ then
\begin{align*}
    P_{\lambda,\Lambda}^+(D^2\varphi(x,t))-\varphi_t(x,t)\geq-\phi(|D \varphi(x,t)|).
\end{align*}
\end{dfn}

Note that the inequalities  \eqref{eq1} and  \eqref{eq2} are not homogeneous.  We thus need the following rescaling lemma from \cite[Lemma 2.4]{A} and  \cite[Lemma 4.4]{j}.
\begin{lemma}\label{scaling}
Let $u\in C(Q_r(x_0,t_0))$ be a viscosity supersolution  of \eqref{eq1}(subsolution of \eqref{eq2}) in $Q_r(x_0,t_0)$. There exists a universal constant $L_2 \ge \Lambda_0$ such that if $A\in (0, \infty)$ then for every $r\leq r_A$, where
\begin{align}\label{scal}
    r_A=\frac{A}{L_2(\phi(A)+A)}=\frac{1}{L_2(\eta(A)+1)}
\end{align}
the rescaled function 
\begin{align*}
    u_r(x,t):=\frac{u(rx,r^2t)}{A},
\end{align*}
is a supersolution of \eqref{eq1} (subsolution of \eqref{eq2}) in its domain.
 \end{lemma}

\section{Proofs of the main results}\label{s:proof}

In this section we give the proofs of the four theorems. 

\begin{proof}[\textbf{Proof of Theorem \ref{thm1}}]
Without loss of generality we assume that $n=3$. The case $n>3$ then follows by adding dummy variables. In the following we write $(x,z) \in \R^{2}\times \R$ for a point in $\R^3$. 

 Fix a small $\varepsilon_0>0$ and choose $q=\frac{4}{\varepsilon_0}$.  For $k \ge k_0(\varepsilon_0)$ denote $r =k^{-\varepsilon_0}$. We  first   define  the functions $u_k$ in $\R^2 \setminus B_r  \times \R$  as
 \[
u_k(x,z)=\frac{1}{k}|x|^{-q}  \quad \text{for } \,  (x,z) \in (\R^2 \setminus B_r) \times \R.
\]
By  direct calculation we have 
\begin{align*}
\partial^2_{x_i x_j}u_k(x,z)=\frac{q(q+2)}{k}|x|^{-q-4} x_i  x_j - \frac{q}{k}|x|^{-q-2}\delta_{ij}.
\end{align*}
The matrix $D^2u_k(x,z)$ has eigenvalues $-\frac{q}{k}|x|^{-q-2}$,  $\frac{q(q+1)}{k}|x|^{-q-2}$ and zero. Therefore, by choosing  $\lambda =1$ and $\Lambda=q+1$, we have 
\[P^-_{\lambda,\Lambda}(D^2u_k)=-(q+1)\frac{q}{k}|x|^{-q-2}+\frac{q(q+1)}{k}|x|^{-q-2}=0.\] 
%

We proceed by  defining $u_k$ in $B_r \times \R$ as 
\begin{align*}
u_k(x,z)=\mathcal A-\mathcal B|x|^2+\mathcal C|x|^4+ \mathcal D(z+2)\rho(|x|),
\end{align*}
where $\rho \in C_{c}^{\infty}([0,r))$ is a cut-off function with $\rho\equiv 1$ in $[0,r/2]$,  $\rho' \le 0$ in $[0,r)$ and $|\rho'|\leq \frac{4}{r}$.  Here, the coefficients $\mathcal A, \mathcal B, \mathcal C$ are chosen so that $u_k$ is a twice differentiable function. By a direct calculations this leads to 
\begin{equation} \label{eq:thm1:ABC}
\mathcal A=\frac{8+6q+q^2}{8k}r^{-q},\ \  
\mathcal B=\frac{q}{4k}(q+4)r^{-q-2},\ \text{and} \
\mathcal C =\frac{q}{8k}(q+2)r^{-q-4}.
\end{equation} 
We then choose 
\begin{equation} \label{eq:thm1:D}
\mathcal D=\frac{r^{-q}}{k}.
\end{equation} 
This choice will be clear later in the argument. By direct calculation we  immediately verify that in $B_r \times (-2,2)$ it holds 
\begin{equation} \label{eq:thm1:hessianbound}
|P^-_{\lambda, \Lambda}(D^2u_k)| \le C({\varepsilon_0}) \frac{r^{-q-2}}{k}. 
\end{equation} 

Notice that in $B_{r/2} \times (-2,2)$, using $\rho \equiv 1,$ we have $|Du_k| \ge |\p_z u_k|= \mathcal D.$ Hence, we find that $|Du_k|^{1+\varepsilon_0} \ge k^{\varepsilon_0}\frac{r^{-q-2}}{k}.$ Here, we  used $r=k^{-\varepsilon_0}$ and $q=\frac{4}{\varepsilon_0}$. 
 Consequently, for $k$ large enough, we conclude that by \eqref{eq:thm1:hessianbound} in $B_{r/2} \times (-2,2)$ it holds 
\[
|P^-_{\lambda,\Lambda}(D^2u_k)| \le |Du_k|^{1+\varepsilon_0},
\]
when $k$ is large enough. 

We  are left to find a lower bound for $|Du_k|$ in $(B_r \setminus B_{r/2}) \times (-2,2).$ Trivially, we have 
\begin{align*}
|Du_k|^2 &\ge |D_xu_k|^2 =|x|^2\left(-2 \mathcal B +4\mathcal C|x|^2 +\mathcal D(z+2)\frac{\rho'}{|x|}\right)^2.
\end{align*}
By the choices of $\mathcal B$ and $\mathcal C$ in \eqref{eq:thm1:ABC} it holds 
\[
 -2 \mathcal B +4\mathcal C|x|^2 \le -2 \mathcal B +4\mathcal C r^2 = -\frac{q}{k}r^{-q-2},
\]
in $B_{r} \times (-2,2)$. Moreover, since $\rho' \le 0$ and $\frac{r}{2} \leq |x| \leq r$ in  $(B_r \setminus B_{r/2}) \times (-2,2)$, we obtain
\begin{align*}
|Du_k|^2 &\ge |x|^2\left(-2 \mathcal B +4\mathcal C|x|^2 \right)^2 \ge \frac{q^2}{4} \frac{ r^{-2q-2}}{k^2}
\end{align*} 
in $(B_r \setminus B_{r/2}) \times (-2,2)$.  Finally, we recall that $q=\frac{4}{\varepsilon_0}$ and $r=k^{-\varepsilon_0}$ to conclude that  
\[
|Du_k|^{1+\varepsilon_0} \ge k^{2\varepsilon_0+\varepsilon_0^2}\frac{r^{-q-2}}{k} 
\]
in $ (B_r \setminus B_{r/2}) \times (-2,2)$. Therefore, by \eqref{eq:thm1:hessianbound}  we have   for $k$ large enough
\[
|P^-_{\lambda,\Lambda}(D^2u_k)| \le |Du_k|^{1+\varepsilon_0} 
\]
in $(B_r \setminus B_{r/2}) \times (-2,2)$ as well.
Thus, we have established that, in $B_2$, 
\[
|P^-_{\lambda,\Lambda}(D^2u_k)| \le |Du_k|^{1+\varepsilon_0}. 
\]

Finally we  notice that  $\inf_{B_1} u_k =\frac{1}{k}$ and $\sup_{B_1} u_k >1.$ Thus 
\[
\int_{\inf_{B_1}u_k}^{\sup_{B_1}u_k} \frac{ds}{s^{1+\varepsilon_0}+s} \ge \frac12\int_{1/k}^{1} \frac{ds}{s} = \frac12 \log(k) \to \infty
\]
as $k \to \infty$. This completes the proof.
\end{proof}
 
\begin{proof}[\textbf{Proof of Theorem \ref{fullharnack}}]
Without loss of generality we may assume that $u>0$. (Otherwise, we define  $u+\varepsilon$ and let $\varepsilon \to 0$.) Let $(x_0,t_0) \in Q_2$ with $Q_{2\rho_0} \cup A^+_{\rho_0}(x_0,t_0) \subset Q_2$ and fix $\rho_0 \le \alpha_0,$ where 
\begin{equation}\label{eq:thm2-1}
\alpha_0=\frac{a_0}{L_2(\eta(\frac{1}{2C})+1)}.
\end{equation}
 Here $A^+_{\rho_0}(x_0,t_0) $ is defined in \eqref{def:Aplus},  $a_0=\frac{1}{C(\eta(u(x_0,t_0))+1)}$ and $C,L_2>1$ are from \eqref{Harnack} and  from Lemma \ref{scaling} respectively. 
Since $u$ is continuous there exists $\rho>0$ such that 
\[
u(x_0,t_0)< 2C \inf_{A^+_{a_0\rho}(x_0,t_0)} u(x,t),
\]
where $C$ and $a_0$ are as above. In order to shorten the notation we denote  $A^+_{a_0\rho}(x_0,t_0)$ by $A^+_{a_0\rho}$. 
We choose the  smallest $\rho$, denoted by $\rho_s$, such that   $a_0\rho_s\le \rho_0$ and  
\begin{equation}\label{cont1}
u(x_0,t_0)= 2C \inf_{A^+_{a_0\rho_s}} u(x,t).
\end{equation}
If such a $\rho_s$ does not exist or if   $a_0\rho_s = \rho_0$, then   the claim of the theorem is trivially true. Let us then assume that  $a_0\rho_s<\rho_0$. From \eqref{cont1} we deduce that  there exists a point $(x_s,t_s)$  in the closure of $A^+_{a_0\rho_s}$ such that  
\[
u(x_0,t_0)= 2C u(x_s,t_s).
\]
We  define 
\[
v(x,t)=2Cu(\rho_sx+x_s,\rho_s^2t+t_s).
\]
Next we recall the definition of $\alpha_0$ in \eqref{eq:thm2-1}. Then $a_0 \rho_s \le \rho_0\leq \alpha_0 $ implies  $\rho_s \le\frac{1}{L_2(\eta(\frac{1}{2C})+1)}.$ Therefore  Lemma \ref{scaling} implies that  $v$ is a positive viscosity supersolution of \eqref{eq1} and viscosity subsolution of \eqref{eq2} in $Q_2$ with $v(0,0)=  u(x_0,t_0).$ Thus, we may  apply \eqref{Harnack}  to  find
\begin{align*}
\sup_{A^-_1} v(a_0x,a_0^2t) \le C v(0,0).
\end{align*}
Recall that $A_1^-=\Big\{(x,t):|x|_{\infty} \leq \frac{c_n}{2}, -1+ \frac{c_n^2}{4} \leq t \leq -1+\frac{c_n^2}{2}\Big\}.$ It is easy to see that $\left(\frac{x_0-x_s}{a_0\rho_s},\frac{t_0-t_s}{(a_0\rho_s)^2}\right) \in A^-_1$. Thus, in particular,  we  have 
\[
2Cu(x_0,t_0)\le C u(x_0,t_0),
\]
 which is a contradiction as $u(x_0,t_0)>0$. Hence, we have $a_0\rho_s=\rho_0$ and 
\[
u(x_0,t_0)\le  2C \inf_{A^+_{\rho_0}} u(x,t).
\]
This completes the proof of the theorem.
\end{proof}

Let us then prove Theorem  \ref{Th:holder}. As usual we use the above Harnack  estimate to deduce that the oscillation has an algebraic decay as we reduce the size of  the cubes. The difficulty is that we need to match the  decay estimate to the size of the associated cubes in order to apply the  intrinsic Harnack inequality. 

\begin{proof}[\textbf{Proof of Theorem \ref{Th:holder}}]
Let $C$ be from Theorem \ref{fullharnack}. We let  $\delta=1-\frac{1}{4C},$ $\rho_0=1$ and  $\omega_0:=\sup_{Q_{1}}u-\inf_{Q_{1}}u$. Note that since $C>1$ then $\delta >\frac12$. We also recall the notation  $r_A:= \frac{1}{L_2(\eta(A)+1)}$ for a number $A>0$ from Lemma \ref{scaling}. For  $k\ge 1$ we  define 
\[
\rho_k =\left(\frac{r_{\delta} c_n}{4C\eta(4)}\right)^k r_{\omega_0/4}.
\]
We claim that for every $k \geq 0$ it holds
\begin{equation}\label{oscde}
\osc_{Q_{\rho_k}}u \le \delta^k \omega_0,
\end{equation}
where $\osc_{Q_{\rho}}u:=\sup_{Q_{\rho}}u-\inf_{Q_{\rho}}u$.

We argue by contradiction. Note that \eqref{oscde} is trivially true for $k=0$ and therefore the contradiction assumption implies that there is $k \geq 0$ such that 
\begin{equation}\label{cl}
\osc_{Q_{\rho_k}}u \le \delta^k \omega_0
\end{equation}
and 
\begin{equation}\label{cg}
\osc_{Q_{\rho_{k+1}}}u > \delta^{k+1} \omega_0.
\end{equation}
We denote $M_k=\sup_{Q_{\rho_k}}u$, $m_k=\inf_{Q_{\rho_k}}u$ and define also  $t_k= \left(1-\frac{c_n^2}{4}\right)\left(\frac{r_{\delta} \rho_k}{2C\eta(4)}\right)^2$. We claim that   one of the following holds
\begin{itemize}
\item[(i)]$M_k-u(0,-t_k) \ge \frac{\delta^k \omega_0}{4}$
\item[(ii)]$u(0,-t_k)-m_k \ge \frac{\delta^k \omega_0}{4}.$
\end{itemize}
Indeed, if both are false then $\osc_{Q_{\rho_k}} u\le \delta^k \omega_0/2.$ Note that from \eqref{cg}, we have $\delta^{k+1} \omega_0 < \osc_{Q_{\rho_k}}u.$ Hence, we obtain $\delta<1/2,$ which is a contradiction. Hence either (i) or (ii) holds.  

Assume first that the condition (i) above holds. We define a function $v :  Q_2 \to \R$ as
\begin{equation}\label{def:thm3-1}
v(x,t)=\frac{M_k-u(r_\delta \rho_k x,(r_\delta \rho_k)^2t -t_k)}{\delta^k\omega_0/4}.
\end{equation}
Let us denote $b_0:=\frac{1}{2C\eta(4)}$. We note that by the choices of $\rho_k, t_k$ and $b_0$ it holds 
\begin{equation}\label{def:thm3-2}
A^+_{b_0r_\delta\rho_k}(0,-t_k) =  Q_{\rho_{k+1}} \subset Q_2, 
\end{equation}
where $A^+_{\rho}$ is defined in \eqref{def:Aplus}.  Let us show that $v$ is a supersolution of \eqref{eq1} and subsolution of \eqref{eq2} in $Q_2 $.

Indeed, by the definition of $\rho_k$ it holds 
\[
r_\delta \rho_k \leq r_\delta^{k+1} r_{\omega_0/4} = \frac{1}{L_2^{k+2}}\frac{1}{(\eta(\delta)+1)^{k+1}}  \frac{1}{\eta(\omega_0/4)+1}
\]
On the other hand, by the condition (P3) in \eqref{eq:P3} on $\eta$ it holds 
\[
\eta(\delta^k \omega_0/4) +1 \leq \Lambda_0^{k+1}   (\eta(\delta) +1)^k (\eta(\omega_0/4) +1). 
\] 
Therefore, since $\Lambda_0 \leq L_2$ it holds 
\[
r_{\delta}\rho_k \leq \frac{1}{L_2(\eta(\delta^k\omega_0/4)+1)}=r_{\delta^k \omega_0/4}
\]
and by Lemma \ref{scaling}, $v$ defined in \eqref{def:thm3-1}   is a supersolution of \eqref{eq1} and subsolution of \eqref{eq2}.

Since $M_k=\sup_{Q_{\rho_k}}u$,  we conclude  that $v$ is non-negative in $Q_2$, by possibly enlarging the constant $L_2$ if necessary. Moreover, the condition (i) and \eqref{cl} imply  $1\le v(0,0) \le 4.$ Therefore since $\eta$ is increasing in $[1,\infty)$ and $\eta \ge 1$, we have  $b_0=\frac{1}{2C\eta(4)} \le \alpha_0 = \frac{1}{C(\eta(v(0,0))+1)}.$
Therefore by Theorem \ref{fullharnack}, we have 
\begin{align*}
\frac{1}{C}v(0,0) \le \inf_{A^+_{b_0}} v(x,t).
\end{align*}
Using $v(0,0)\geq 1$, we have  
\begin{align*}
\frac{\delta^k\omega_0}{4C} \le M_k -\sup_{A^+_{b_0}} u(r_\delta \rho_k x,(r_\delta \rho_k)^2t -t_k).
\end{align*}
Recalling \eqref{def:thm3-2}  the above yields 
\begin{align*}
\sup_{Q_{\rho_{k+1}}}u-m_k \le M_k -m_k - \frac{\delta^k\omega_0}{4C} = \osc_{Q_{\rho_{k}}}u - \frac{\delta^k\omega_0}{4C} \le \delta^k \omega_0 - \frac{\delta^k\omega_0}{4C}=\delta^{k+1} \omega_0.
\end{align*}
Hence, we obtain 
\begin{align*}
\osc_{Q_{\rho_{k+1}}}u \le \delta^{k+1} \omega_0
\end{align*} 
which  contradicts  \eqref{cg}. 

If  the condition (ii)  holds,  we define 
\[
v(x,t)=\frac{u(r_\delta \rho_kx,(r_\delta \rho_k)^2t -t_k)-m_k}{\delta^k\omega_0/4}.
\]
By repeating the above argument leads again to a  contradiction. (We leave the details for the reader). Therefore we finally have the estimate  \eqref{oscde}.

Let us now prove the H\"older continuity.  Given $r\le 1$,  choose $k$ such that  
\begin{align*}
\left(\frac{r_{\delta} c_n}{4C\eta(4)}\right)^{k+1} \le r \le \left(\frac{r_{\delta} c_n}{4C\eta(4)}\right)^k
\end{align*}
and use \eqref{oscde} to conclude 
\begin{align*}
\osc_{Q_{r}}u\left(r_{\omega_0/4}x,r_{\omega_0/4}^2t\right) \le \frac{1}{\delta}r^{\alpha} \omega_0,
\end{align*}
where $\alpha =\min\left\{\frac{1}{2}, \frac{\log \delta}{\log \left(\frac{r_{\delta} c_n}{4C\eta(4)}\right)}\right\}$ is a universal constant. Therefore, we have  for all $(x,t) \in Q_1$ 
\begin{align*}
\big|u\left(r_{\omega_0/4}x,r_{\omega_0/4}^2t\right)-u(0,0)\big| \le \frac{1}{\delta}(|x|+|t|^{1/2})^{\alpha} \omega_0.
\end{align*}
Thus, for $(x,t) \in Q_{r_{\omega_0/4}}$ we have 
\begin{align*}
|u\left(x,t\right)-u(0,0)| \le \frac{1}{r_{\omega_0/4}^{\alpha}\delta}(|x|+|t|^{1/2})^{\alpha} \omega_0 \le \tilde C \phi(\omega_0)(|x|+|t|^{1/2})^{\alpha},
\end{align*}
where in the second inequality we have used the fact that $\phi(s) = \eta(s) s$ for $\eta \geq 1$  and  $\alpha <1$. Here $\tilde C>1$ is a universal constant.  On the other hand, it is easy to see that for all $(x,t) \in Q_1 \setminus Q_{r_{\omega_0/4}}$ it holds 
\begin{align*}
|u\left(x,t\right)-u(0,0)|  \le \omega_0 \frac{r_{\omega_0/4}^{\alpha}}{r_{\omega_0/4}^{\alpha}} \le  \tilde C \phi(\omega_0)(|x|+|t|^{1/2})^{\alpha}.
\end{align*}
Therefore, we conclude that for all  $(x,t) \in Q_{1},$ it holds
 \begin{align*}
|u\left(x,t\right)-u(0,0)|  \le\tilde  C \phi(||u||_{\infty})(|x|+|t|^{1/2})^{\alpha}.
\end{align*}
 This implies the claim by using the fact that the equations \eqref{eq1} and \eqref{eq2} are translation invariant together with a standard covering argument.

\end{proof}
\begin{proof}[\textbf{Proof of Theorem \ref{Th:int_h}}]
Without loss of generality we may assume that $u>0$. (Otherwise, we define  $u+\varepsilon$ and let $\varepsilon \to 0$.) 

We define a sequence of  radii $r_i>0$  and time levels $t_i<0$ as  follows.  Set first $M_0=u(0,0)$, $\rho_1=\frac{1}{C(\eta(M_0)+1)}$ and choose the radius  $r_1=\frac{\rho_1c_n}{2}$ and the time  $t_1=\left(-1+\frac{c_n^2}{4}\right)\rho_1^2$. Here $C>1$ and $c_n \leq 1$ are from \eqref{Harnack}.  We proceed by defining $r_i$ iteratively such that if $r_{i}$ is defined, we let 
\begin{equation}\label{def:M-i}
\tilde M_i= \sup_{x \in B_{r_i}}u(x,t_i) \quad \text{and} \quad M_i = \max \{ M_{i-1}, \tilde M_i\},
\end{equation}
  set $\rho_{i+1}=\frac{1}{C(\eta(M_i)+1)}$ and define 
\[
r_{i+1}=r_i+\frac{\rho_{i+1}c_n}{2} \quad \text{and} \quad t_{i+1}=t_i+\left(-1+\frac{c_n^2}{4}\right)\rho_{i+1}^2.
\]

We first observe that there exists $K \in \mathbb{N}$ such that $r_{K} \ge 1$. Indeed, if $r_i < 1$ for all $i \in \mathbb{N}$ then the definition of $r_i$ yields
\[
\frac{c_n}{2C} \sum_{i=0}^{\infty} \frac{1}{\eta(M_i)+1} \leq 1 
\]
and $t_i \geq -1$ for all $i$. But then necessarily $\rho_i \to 0$ as $i \to \infty$, which means that $\eta(M_i) \to \infty$ as $i \to \infty$. This in turn implies $M_i \to \infty$, which is a contradiction since $u$ is continuous. We define $K$ to be the first index for which $r_K \geq 1$ and note that $r_K \leq 2$.

We proceed by claiming that for all $i = 0,1, \dots, K-1 $ it holds 
\begin{equation}\label{thm4-1} 
M_{i+1} \le CM_i.
\end{equation}
If $M_{i+1} = M_i$ then \eqref{thm4-1}  is trivially true. Let us then assume that  $M_{i+1}> M_i$.  We  choose $x_{i+1} \in \overline{B}_{r_{i+1}}$ such that   $M_{i+1}=u(x_{i+1},t_{i+1})$ and let $y_{i}$ be the closest point to $x_{i+1}$ in $\overline{B}_{r_i}$. Then it holds $u(y_i, t_i) \leq M_i$.  We define 
\[
f(r):=u\left(x_{i+1}+\frac{r}{\rho_{i+1}}(y_{i}-x_{i+1}),{t}_{i+1}+\left(1- \frac{c_n^2}{4}\right)r^2\right)
\]
and notice that $f(0)= u(x_{i+1},t_{i+1}) =M_{i+1} >M_{i}\geq u(y_i,t_i) = f(\rho_{i+1})$. By continuity  of $u$ there exists $0<\rho_0 \le \rho_{i+1}$ such that $f(\rho_0)=M_i$. We set 
\[
(x_0,t_0)=\left(x_{i+1}+\frac{\rho_0}{\rho_{i+1}}(y_{i}-x_{i+1}),{t}_{i+1}-\left(-1+\frac{c_n^2}{4}\right)\rho_0^2\right)
\]
 and apply \eqref{Harnack} to have 
\[
 \sup_{A_{\rho_0}^-(x_0,t_0)}u(x,t) \le CM_i.
\]
 Notice that $(x_{i+1},t_{i+1}) \in A_{\rho_0}^-(x_0,t_0)$ and therefore it holds $ u(x_{i+1},t_{i+1}) =  M_{i+1} \le CM_i$, which proves \eqref{thm4-1}.

We use \eqref{thm4-1}  to estimate
\[
	\int_{M_0}^{M_{K}} \frac{dt}{\phi(t)+t} 
 	\le \sum_{i=0}^{K-1}\int_{M_i}^{M_{i+1}} \frac{dt}{\phi(t)+t}
 	\le \sum_{i=0}^{K-1}\int_{M_i}^{CM_i} \frac{dt}{\phi(t)+t}.
 \]
Since  $\phi$ is increasing and of the form $\phi(s) = \eta(s)s$ we obtain 
\[
\begin{split}
 	\sum_{i=0}^{K-1}\int_{M_i}^{CM_i} \frac{dt}{\phi(t)+t} &\le \sum_{i=0}^{K-1}\int_{M_i}^{CM_i} \frac{dt}{\phi(M_i)+M_i}
 	\\
&=\sum_{i=0}^{K-1} \frac{(C-1)M_i}{\phi(M_i)+M_i } =(C - 1)\sum_{i=0}^{K-1} \frac{1}{\eta(M_i)+1 } \\
&=  C(C-1)\sum_{i=0}^{K-1}\rho_{i+1} \leq \frac{2C^2}{c_n} r_K.
\end{split}
 \]
Since $r_K \leq 2$ and $\phi(t) \geq t$, we conclude 
 \begin{equation}\label{eq:int}
 \int_{M_0}^{M_K} \frac{dt}{\phi(t)} \le \tilde C,
 \end{equation}
for a universal constant $\tilde C$.

Note that in \eqref{eq:int} we have  $M_0 = u(0,0)$  and $M_K \geq \sup_{x \in B_1}u(x,t_K)$. We let $t_2 = -t_K$ and need yet to estimate $|t_K|$. 
By the definition of $t_K$ and $\rho_i$ it holds 
 \begin{equation}\label{eq:bound-t-k}
|t_K| = \left(1 - \frac{c_n^2}{4} \right) \sum_{i=1}^{K}\rho_i^2 \leq  \sum_{i=1}^{K}\rho_i^2. 
 \end{equation}
Then by increasing $C$  if necessary such that $C^{-1}\leq \frac{c_n}{4}$,  we have $\rho_i = \frac{1}{C(\eta(M_i)+1)} \leq \frac{c_n}{4}$. Then since 
$r_K \leq 2$ we have
\[
|t_K| \leq  \sum_{i=1}^{K}\rho_i^2 \leq \frac{c_n}{4} \sum_{i=1}^{K}\rho_i = \frac{r_K}{2} \leq 1.
\]
On the other hand, trivially it  holds
\[
|t_K| \geq |t_1| \geq \frac{1}{2C^2(\eta(u(0,0))+1)^2}. 
\]
 Consequently, we obtain 
\[
\frac{1}{\tilde C(\eta(u(0,0))+1)^2} \le t_2 \leq 1,
\]
for a universal constant $\tilde C$. 

 Similarly, we  prove 
 \[
 \int^{u(0,0)}_{\inf_{B_1}u(\cdot,t_1)} \frac{dt}{\phi(t)} \leq \tilde C
\]
for a universal constant  $\tilde C$ and a time level $t_1>0$ with  
\[
\frac{1}{\tilde C (\eta(u(0,0))+1)^2} \leq t_1 \leq 1.
\]
 This completes the proof.
\end{proof}

Corollary  \ref{cor:harnack} follows from Theorem \ref{Th:int_h} but we include the argument for the reader's convenience. 

\begin{proof}[\textbf{Proof of  Corollary \ref{cor:harnack} }]
Let $u \in C(Q_4(0,1))$ be as in the statement and define $u_{\varepsilon} = u+ \varepsilon$.  We adopt the notation from the proof of Theorem \ref{Th:int_h} and denote $M_i(\varepsilon)$ as  in \eqref{def:M-i} for $u_{\varepsilon}$ and $\rho_{i+1} = \frac{1}{C(\eta(M_i(\varepsilon))+1)}$. Since $u_{\varepsilon}(0,0) = \varepsilon$ we have by Theorem \ref{Th:int_h}  that 
 \[
 \int_{\varepsilon}^{M_K(\varepsilon)} \frac{dt}{\phi(t)} \leq  C
\]
for $M_K(\varepsilon) \geq \sup_{x \in B_1}u_\varepsilon(x,t_K(\varepsilon))$. Here $t_K(\varepsilon) = -t_2$.  But since $\phi$ satisfies \eqref{eq:osgood} then $M_K(\varepsilon) \to  0$ as $\varepsilon \to 0$. We need yet to prove that we may choose $t_K(\varepsilon)$ such that  $t_K(\varepsilon) \to  0$ as $\varepsilon \to 0$. 

We divide the proof in two cases. Assume first
\[
\lim_{s \to 0} \eta(s) < \infty. 
\]  
Then there is $\tilde C$ such that $1 \leq \eta(s) \leq \tilde C$ for all $s \in (0,1]$. Therefore the equations \eqref{eq1} and  \eqref{eq2} reduce to the linear case and  we have the classical Harnack inequality by \cite{is, W1} and the claim follows.  

Assume then 
 \[
\lim_{s \to 0} \eta(s) = \infty. 
\]  
 The definition \eqref{def:M-i} of $M_i(\varepsilon)$ implies  $M_i(\varepsilon) \leq M_K(\varepsilon)$. Since  $M_K(\varepsilon) \to  0$ as $\varepsilon \to 0$, then for every $\delta>0$ it holds  $\rho_{i+1} = \frac{1}{C(\eta(M_i(\varepsilon)) +1)} <\delta$ for all $i \leq K(\varepsilon)-1$ when $\varepsilon$ is small. By the definition of $r_{K(\varepsilon)}$ it holds 
\[
\sum_{i=1}^{K(\varepsilon)}\rho_i  = \frac{2}{c_n} r_{K(\varepsilon)} \leq  \frac{4}{c_n}.
\]
Therefore by \eqref{eq:bound-t-k} we have 
\[
|t_K(\varepsilon)| \leq  \sum_{i=1}^{K(\varepsilon)}\rho_i^2 \leq \delta  \sum_{i=1}^{K(\varepsilon)}\rho_i  \leq  \frac{4\delta}{c_n} .
\]
The claim follows by letting first $\varepsilon \to 0$ and then $\delta \to 0$ .
\end{proof}

\section*{Acknowledgments}
The authors were supported by the Academy of Finland grant 314227.

\end{document}